
\magnification=\magstep1
\voffset=1.00truein
\settabs 18 \columns
\hoffset=0.00truein
\baselineskip=17 pt

\font\stepone=cmr10 scaled\magstep1
\font\steptwo=cmr10 scaled\magstep2

\def\b{\bigskip}
\def\bb{\bigskip\bigskip}
\def\bbb{\bigskip\bigskip\bigskip}

\def\no{\noindent}

\def\ce{\centerline}
\def\ve{\vfill\eject}

\def\e e{$e^+ e^-$ }

\def\A{\hbox{$A(\langle q\rangle,a)$}}

 \def\q{\quad}

\bbb
\ce{\steptwo {Deformations of Multiparameter Quantum $gl(N)$}}
\bb
\ce {\bf {C. Fr{\o}nsdal}}
\b
\ce {Department of Physics, UCLA, Los Angeles, California 90024-1547}
\bb
\ce {\bf {A. Galindo}}
\b
\ce {Departamento de Fisica Te\'{o}rica}
\ce {Universidad Complutense, 28040 Madrid, Spain}
\bbb
\ce {\bf Abstract}
\b
\no Multiparameter quantum $gl(N)$ is not a rigid structure.  This paper
defines
an essential deformation as one that cannot be interpreted in terms of a
similarity transformation, nor as a perturbation of the parameters.  All the
equivalence classes of first order essential deformations are found, as well as
a class of exact deformations.  This work provides quantization of all the
classical Lie bialgebra structures (constant r-matrices) found by Belavin and
Drinfeld for $sl(n)$.  A special case, that requires the Hecke  parameter to be
a cubic root of unity, stands out.
\vfill\eject

\no {\stepone 1. Introduction}

Belavin and Drinfeld [1] classified the r-matrices (structures of coboundary
Lie bialgebra) associated to simple Lie algebras, both finite and affine.   A
program of ``quantization" of Lie algebras proposed by Drinfeld [2] would
promote the classical structures to bialgebra deformations of enveloping
algebras.  Standard forms of such ``quantum deformations" of the simple affine
Lie algebras were obtained by Jimbo [3]. Here we are interested only in the
finite case, of constant r- and R-matrices.  So far there is no general
classification of quantized  Lie algebras.

Deformation theory applied at the classical point [4] is difficult, since the
obstructions appear only in the second order.  But a large family of exact
quantum deformations of $gl(N)$ is known [5], with $1 + N(N-1)/2$ parameters.
These algebras are rigid (with respect to essential deformations) at generic
points in parameter space, even  to first order deformations, but essential
deformations exist on algebraic surfaces of lower dimension, for $N>2$.  The
determination of all the first order deformations, presented here, goes far
towards a complete classification of all formal and/or exact deformations.  All
first order deformations are combinations of ``elementary" deformations, and
all
the elementary deformations are exact.  A family of exact deformations that
contains the largest possible number of elementary ones is described.  It is
characterized by the fact that the dual Lie algebra has the largest possible
semi-simple radical, which is the main reason why these new quantum groups are
of interest in connection with integrable models of quantum field theories and
classical statistical systems.  We have in mind applications to new Toda field
theories, statistical models and knot theory.

\q {\bf Added Notes.} (1) This paper was distributed as a  UCLA preprint in
December, 1992.  We have since become aware of a paper by Cremmer and Gervais
[11]   in which they describe an R-matrix that is substantially the same as the
one  we present in Section 9. (2) The Hecke  condition that was imposed on the
deformations considered in the  original version of this paper is not
essential; if  the  Hecke parameter $a \neq \pm 1$,  then all the deformations
that  preserve the Yang-Baxter relation provide representations of the braid
group  factored through the Hecke algebra. The text has been modified to
include
this observation. (3) The cohomology that underlies the  structure of the
deformation space is not explained in this paper.  However, the corresponding
problem for Lie bialgebras has been clarified in [12]. We hope to return to the
question in a future  publication.

\bbb
\no {\stepone 2. Multiparameter Quantum $gl(N)$ and Quantum P-Algebras}
\b We consider the free associative algebra ${\cal F}_x$ generated by $(x^i),
\,\, i =
\cdots, N$, and the ideal ${\cal F}_{x0}$ generated by
$$x^ix^j - q^{ij}x^jx^i, \,\, i,j=1,\cdots ,N ,\eqno(2.1)
$$
\no in which the $q$'s are taken from a field K with characteristic 0, with
$$ q^{ij}q^{ji} = 1, \,\, q^{ii} = 1, \,\, i,j = 1, \cdots ,N . \eqno(2.2)
$$
\no We call {\bf quantum plane} the associative algebra
${\cal F}_x/ {\cal   F}_{x0}$; that is, the associative algebra generated by
the
$x$'s with relations
$$ x^ix^j - q^{ij}x^jx^i = 0, \,\, i,j = 1, \cdots ,N . \eqno(2.3)
$$
\no Similarly, the quantum anti-plane
${\cal F}_{\theta}/ {\cal F}_{\theta 0}$ is an associative algebra generated by
N elements
$(\theta^i), \,\, i = 1, \cdots , N$, with relations
$$
\theta^i \theta^j + r^{ij} \theta^j \theta^i = 0, \,\, i,j = 1, \cdots ,N ,
\eqno(2.4)
$$
\no in which the $r$'s are parameters from $K$  satisfying the same relations
as
the $q$'s, Eq. (2.2).  Let $V$ denote the linear vector space over $K$ spanned
by the $x$'s (or by the $\theta$'s).
\b
\no {\bf Definition}.  A {\bf generalized symmetry} is an element $P$ of
${\rm End}(V \otimes V)$ that satisfies the Hecke condition
$$ (P-1)(P+a) = 0, \eqno(2.5)
$$
\no for some $a \in K, a \not= -1,0$.  We shall say that the tensor $xx =
(x^ix^
j)$ is
$P$-symmetric, and that the tensor $\theta \theta = (\theta^i \theta^j)$ is
$P$-antisymmetric, if
$$xx(P-1) = 0, \,\,\, \theta \theta(P+a) = 0 . \eqno(\hbox{2.3'-4'})$$
\b Let $P_{12}$ be the operator on $V \otimes V \otimes V$ that acts as $P$ on
the two first factors; the braid relation is
$$ P_{12} P_{23} P_{12} = P_{23} P_{12} P_{23} . \eqno (2.6)
$$
\b
\no{\bf Theorem 1}. Given relations (2.3) and (2.4), with parameters $q$ and
$r$ subject to the conditions (2.2),  the following two statements are
equivalent:
$${\rm (a)} \,\,\,
   r^{ij} = aq^{ij}, \,\, i < j, \,\, i,j = 1, \cdots, N ; \eqno(2.7)
$$
\no{\rm (b)} There exists $P$ in ${\rm End}(V \otimes V)$  satisfying (2.5) and
(2.6),  such that the relations (2.3-4) coincide respectively with Eqs.
(2.3'-4');  it is unique  up to a permutation of the basis.
\b Let $P$ be a generalized symmetry of dimension $N$.   Consider the algebra
${\cal F}_x$ freely generated  by $(x^i), \,\, i=1, \cdots ,N$, with the ideal
${\cal F}_{x0}$ generated by
$$ (xx(P-1))^{ij}, \,\, i,j = 1, \cdots, N;\eqno(2.8)
$$
\no and the algebra ${\cal F}_{\theta}$ freely generated by
$(\theta^i) \,\, i= 1, \cdots , N$,  with the ideal ${\cal F}_{\theta 0}$
generated by
$$
\theta \theta (P+a).\eqno(2.9)
$$
\no Let ${\cal F}$ be the algebra generated by
$(x^i)$ and $(\theta^i), \,\, i= 1, \cdots , N$ with relations
$$ ax \theta = \theta xP . \eqno(2.10)
$$
\no This algebra contains ${\cal F}_x$ and ${\cal F}_{\theta}$ as subalgebras
and the ideals ${\cal F}_{x0}$ and ${\cal F}_{\theta 0}$ are thus canonically
identified with subsets of ${\cal F}$.
\b
\no {\bf Theorem 2.}  Suppose that the generalized symmetry $P$ satisfies the
braid relation.  Let $X$ denote the linear span of the $x$'s and $\Theta$ the
linear span of the $\theta$'s, then the statements
$$
\eqalignno{ {\cal F}_{x0} \Theta = \Theta {\cal F}_{x0}&&(2.11)\cr {\cal
F}_{\theta 0} X = X {\cal F}_{\theta 0}&&(2.12)\cr}
$$
\no hold in ${\cal F}$.  Conversely, if  both (211) and (2.12) hold, then $P$
satisfies the braid relation.
\b
\no {\bf Proof.}  Eqs. (2.11), (2.12) are equivalent, respectively, to
$$ ({\rm braid})_{123} (P_{12}-1) = 0, \quad  ({\rm braid})_{123} (P_{12}+a) =
0, \eqno(2.13)$$ with
$$  ({\rm braid})_{123}:= P_{12} P_{23} P_{12} - P_{23} P_{12} P_{23}.
\eqno(2.14)
$$
\b
\no {\bf Definition.}  Let $\langle q\rangle$ stand  for a set of parameters
$(q^{ij}), \,\, i,j = 1, \cdots, N$, satisfying $q^{ij}q^{ji} = 1$ and $q^{ii}
=
1$; and $a$ an additional parameter, all in the field $K$.   Let $r^{ij} =
aq^{ij}$ for $i<j, \,\, r^{ii} = 1$ and $r^{ji} = q^{ji}/a$ for $i<j$.  The
{\bf
standard quantum algebra} {\A}  is generated by the $x$'s and the $\theta$'s,
with relations (2.3, 2.4) and (2.10).  More generally,  for any generalized
symmetry $P$, the algebra $A(P)$ generated by the $x$'s and the $\theta$'s,
with
relations (2.3'-4') and (2.10), will be called a {\bf quantum P-algebra}  if
the
conditions (2.11) and (2.12) hold.
\b The original aim of this work was to calculate the  deformations of {\A}  in
the category of quantum P-algebras. Later it was discovered that the Hecke
condition is inessential; it is enough to impose the braid relation,  then the
Hecke condition is satisfied automatically. Therefore we have, in fact,
determined all the (equivalence classes of essential) deformations in the
category of Yang-Baxter matrices. It is interesting that all the deformations
lend  themselves to an interpretation in terms of generalized symmetry.  The
quantum pseudogroup (in the sense of Woronowicz [6]) associated to $P$ is the
unital algebra generated by the matrix elements of an $N$-by-$N$ matrix $T$,
with relations
$$ [P,T\otimes T]=0.\eqno(2.15)
$$
\no It is the algebra of linear automorphisms of $A(P)$; that is, the  set of
mappings
$$ (x,\theta,T)\rightarrow x\otimes T, \,\,\, \theta \otimes T\eqno(2.16)
$$
\no that preserve the relations (2.3'-4', 2.10) of $A(P)$.   It is related, via
duality, to a quantum group in the sense of Drinfeld [2].  Twisted, quantum
$gl(N)$ [5] corresponds to ${\rm Aut}\A$.  The deformations of this quantum
group are thus in 1-to-1 correspondence with the deformations of the standard
quantum algebra \A.

The next section defines the deformations that will be calculated in this
paper.
Then we begin, in Section 4, by calculating the first order deformations.   It
turns out that {\A} is rigid for parameters in general position.  Interesting
nontrivial deformations (even exact ones) do exist on certain algebraic
surfaces
in parameter space, for $N>2$.  The exact nontrivial deformations include the
``esoteric" quantizations of $gl(N)$ reported elsewhere;   they are described
in
Section 6. The existence of an unexpected special case  that requires $a^3 = 1$
deserves some attention.
\ve

\no {\stepone 3. Deformations}
\b Henceforth,  $P$ will denote the generalized symmetry associated with \A. A
{\bf formal deformation} of {\A}   is here a quantum $P(\epsilon)$-algebra
with
$P(\epsilon)$ a formal power series in an indeterminate $\epsilon$
$$ P(\epsilon)=P + \epsilon P_1 + \epsilon^2 P_2 + \cdots ,\eqno(3.1)
$$
\no that satisfies the Hecke condition with the parameter $a$ independent of
$\epsilon$, and such that (2.11, 2.12) hold.  In this case we shall say that
$P(\epsilon)$ is a formal deformation of $P$.  A deformation is {\bf exact} if
the series $P(\epsilon)$ has a nonvanishing radius of convergence.

If $P(\epsilon)$ is a formal deformation of $P$, then
$$ P(\epsilon,1)=P+\epsilon P_1\eqno(3.2)
$$
\no is a first order deformation.  More generally, a {\bf first order
deformation} is defined as a formal deformation except that one sets
$\epsilon^2 = 0$.  A first order deformation is not necessarily the first two
terms of a formal deformation.  For example, at the classical point, where $a$
and all the $q$'s are equal to unity, the braid relation is moot in first
order.  (Recall that the condition that defines a formal deformation is in this
case the classical Yang-Baxter relation, which is second order.)  For this
reason, the concept of a first order deformation is of no use at the classical
point.  In contrast with this, we shall find that, at general position in
parameter space, {\A}   is rigid with respect to first order deformations,
which
implies rigidity under formal deformations.  The case $a=1$, with the
$q$'s in general position, is in this respect intermediate; it is best treated
separately.

Two types of deformations (and combinations of them) will be considered
trivial.  A linear transformation, with coefficients in $K[\epsilon]$,
$$ x^i \rightarrow x^i + \epsilon x^j A_j^i + \cdots , \,\, \theta^i
\rightarrow
\theta^i + \epsilon \theta^j A_j^i + \cdots , \eqno(3.3)
$$
\no induces a trivial, formal deformation of $P$.  A variation of the $q$'s
$$ q^{ij} \rightarrow q^{ij} + \epsilon \delta q^{ij} + \cdots \eqno(3.4)
$$
\no will also be considered trivial.  A deformation that is not trivial is
called essential.  Two first order deformations $P_1$ and $P'_1$ are
equivalent
if the difference is trivial; that is, induced by  transformations of the type
(3.3), (3.4). We shall classify the equivalence  classes of first order
deformations.  The following  result is deduced from an examination of the case
$N=2$ in the next section.
\b
\no {\bf Theorem 3.}  If $a \not= 1$, then each equivalence class of first
order
deformations contains a unique representative with the property that
$(P_1)_{kl}^{ij} = 0$ for every index set $i,j,k,l$ that contains no more than
two different numbers.
\b The first order deformation of $P$ induced by (3.3) is
$$ P_1 = PZ - ZP, \quad Z:= A \otimes 1 + 1 \otimes A,\eqno(3.5)
$$
\no or more explicitly
$$ (P_1)_{kl}^{ij} = a(\hat{q}^{lk} - \hat{q}^{ij}) Z_{kl}^{ji} + (1-a)[(k<l)
 - (i<j)] Z_{kl}^{ij}, \eqno(3.6)
$$
\no where
$$ {
\hat{q}^{ij}:=\cases{ q^{ij}&{\rm if} $i<j$,\cr q^{ij}/a&{\rm if} $i\geq
j$,\cr}
\qquad (i<j):=\cases{ 1&{\rm if} $i<j$,\cr 0&{\rm otherwise}.\cr} }
\eqno(3.7)
$$ Preservation of the Hecke condition (2.5) under first order deformations is
equivalent to requiring that
$$
\eqalign{q^{ij} (P_1)_{lk}^{ji} + \hat{q}^{kl} (P_1)_{kl}^{ij} &=0 \cr (1-a)
q^{ij} (P_1)_{lk}^{ji} =(a-1) \hat{q}^{kl} (P_1)_{kl}^{ij}  &= (P_1)_{lk}^{ij}
+
aq^{ij} \hat{q}^{kl} (P_1)_{kl}^{ji},
 \,\, i \leq j, \,\, k \leq l.\cr}\eqno(3.8)
$$
\no The main difficulty is to extract the conditions on $P_1$ imposed by the
braid relation.  We found that the strategy made available by Theorem 2
simplifies this task.  Both conditions, (2.11) and (2.12), must be invoked.  We
leave out the details.
\ve

\no {\stepone 4. The Case $N=2$}
\b The components of $P_1$ come in three disjoint sets.  We comment first on
the
generic case, $a$ and $q:=q^{12}$ in general position.

\no (1) Since $a\not= -1$, the Hecke condition requires
$$ (P_1)_{11}^{11}=(P_1)_{22}^{11}=(P_1)_{11}^{22}=(P_1)_{22}^{22}=0.\eqno(4.1)
$$
\no (2) The components $(P_1)_{12}^{12}, \, (P_1)_{12}^{21}, \,
(P_1)_{21}^{12}, \, (P_1)_{21}^{21}$ are unaffected by the transformations
(3.3).   The implications of the Hecke condition and the braid relation (to
order $\epsilon$) are
$$ (P_1)_{12}^{12}+(P_1)_{21}^{21}=0,\,\,  (a-1)(P_1)_{21}^{21} = 0, \,\,
(P_1)_{21}^{12} + aq^2(P_1)_{12}^{21} = 0.\eqno(4.2)
$$
\no This allows for a one-dimensional space of deformations that amounts to a
variation of $q$; hence it is trivial.

\no (3) The remaining coefficients are affected by the transformations (3.3).
There are eight of them, of which 4 remain independent after imposition of the
Hecke condition.  Eq. (2.11) imposes two additional constraints, namely
$$
\eqalign{(q-1)(P_1)_{12}^{11} &+ (aq-1) (P_1)_{22}^{21} = 0, \cr
           (q-1) (P_1)_{21}^{22} &+ (q-1/a) (P_1)_{11}^{12} = 0.\cr}
\eqno (4.3)
$$
\no This leaves two free parameters; the deformations are induced by
transformations of the type (3.3), so they are trivial.  No additional
conditions from (2.12).

The exceptional cases are as follows.  If $a=q=1$ we get nothing at all from
the
braid relation, to first order.  The classical Yang-Baxter relation appears in
the second order.  The only other special case is $a=1, \,\, q\not=1$.   We now
have two additional free parameters, say $(P_1)_{21}^{12}$ and
$(P_1)_{12}^{12}$.  The former is interpreted as a variation of $q$; the latter
is nontrivial.  No additional restrictions from imposing (2.12).

The space of nontrivial, first order deformations is thus one-dimensional in
the
exceptional case when $a=1, \, q \not= 1$.  But this first order deformation is
not the first order part of a formal (power series) deformation.  {\bf The
space
of essential, formal deformations of {\A}  is empty in the case
$N=2$, unless $a = q = 1$.}

\ve

\no {\stepone 5. The Case $N=3$}
\b We shall assume, from here onwards, that $a^2 \not= 1$. In view of  Theorem
3, we may ``fix the gauge" by setting
$(P_1)_{kl}^{ij}=0$ unless all three of the values $1,2,3$ appear among the
indices.

We find that there are two very different cases:

(i) When the repeated index $k$ lies between the other two, $i<k<j$, then the
complete set of conditions is
$$
\matrix{ (P_1)_{ij}^{kk}=-aq^{ji}(P_1)_{ji}^{kk}\hfill&\neq 0
\Rightarrow& q^{ik}=q^{kj}\quad{\rm and}\quad q^{ij}=(q^{ij})^2,\hfill\cr
(P_1)_{kk}^{ij}=-q^{ij}(P_1)_{kk}^{ji}\hfill&\neq 0
\Rightarrow& q^{ik} = q^{kj}\quad {\rm and}\quad q^{ij}=a(q^{ik})^2.
\hfill\cr}
\eqno (5.1)
$$ (ii) In the other cases one finds:
$$
\matrix{ (P_1)_{ij}^{kk}=-aq^{ji}(P_1)_{ji}^{kk}\hfill&\neq 0
\Rightarrow&\cases{ q^{ik}=a^2q^{kj},q^{ij}=(q^{kj})^2& {\rm if} $k<i<j$,\cr
q^{ik} =a^{-2}q^{kj},q^{ij} = (q^{ik})^2& {\rm if} $i<j<k$,\cr}   }\eqno (5.2)
$$
$$
\matrix{ (P_1)_{kk}^{ij}=-q^{ij}(P_1)_{kk}^{ji}\hfill&\neq 0
\Rightarrow&\cases{ q^{ik}=q^{kj},q^{ij}=a(q^{kj})^2& {\rm if} $k<i<j$,\cr
q^{ik} = q^{kj},q^{ij} =a(q^{ik})^2& {\rm if} $i<j<k$,\cr}   }\eqno (5.3)
$$
\no and, in addition, very unexpectedly, in this case
$$ a^3 = 1. \eqno (5.4)
$$
\b
\no {\stepone 6. The Case $N=4$}
\b The deformations found for $N=3$ give rise to a class of deformations for
general $N$ that may be called class-3 deformations.  [The negative result for
$N=2$ means that there are no class-2 deformations.]  These class-3
deformations
will be discussed later; here we investigate those of class 4; that is, those
that involve $(P_1)_{kl}^{ij}$'s with four different index values.  They are
all
essential (first order) deformations.

The Hecke condition and the braid relation connect $P_1$'s with the same set
$(i,j)$ of upper indices and the same set $(k,l)$ of lower indices to each
other.  There are six sets, distinguished by the relative order of the four
indices.  We fix $i<j$ and $k<l$.  The result is that deformations satisfying
(2.11), (2.12) exist in 2 of the 6 cases, the conditions being,  for $a^2 \not=
1$:
$$
\matrix {k<i<j<l:& x=y=1/a, &u=v=1,&\delta P_{kl}^{ij} = 0,\cr
 i<k<l<j:& x=y=a, &u=v=1,&\delta P_{kl}^{ij} = 0,\cr}\eqno(6.1)
$$
\no with
$$
\matrix{ x:= q^{ij} q^{jk} q^{jl},y:=q^{ij} q^{ki} q^{li}, u:=q^{kl} q^{lj}
q^{li},v:=q^{kl} q^{ik} q^{jk},xu = yv .\cr} \eqno(6.2)
$$
\no The space of class-4 deformations is thus at most 1-dimensional (for
$N=4$).
\bbb
\no {\stepone 7. First Order Deformations in the General Case}
\b When $N>3$ there are class-3 deformations involving $(P_1)_{ij}^{kk}$ and/or
$(P_1)_{kk}^{ij}$, for each triple $i,j,k$ of unequal numbers in $1, \cdots,
N$.  But there are now additional conditions that must be satisfied by the
additional $q$'s, and $i,j,k$ must be nearest neighbours.

When $N>4$ there are class-4 deformations satisfying (6.1).  But the braid
relation imposes additional conditions on the additional $q$'s, and the
position
of $i,j,k,l$ inside the set $1,2 \ldots N$ is restricted so that either $k+1 =
i<j = l-1$, or else $i+1 = k<l = j-1$.
\b
\no {\bf Definition.}  An elementary first order deformation is one in which
some $(P_1)_{kl}^{ij}$ is non-zero for just one unordered pair $i,j$  and just
one unordered pair $k,l$.
\b
\no {\bf Theorem 4.}  Suppose $a \neq \pm 1, 0 $. There are two series of
elemen
tary, first order, essential deformations. The ``principal series"  is
described
first. Let $i,j$ be any index pair with, either (case 1)
$ k + 1 = i \leq j = l - 1$, or else (case 2) $i+1 = k \leq l = j-1$.  Let
$P_1=0$, except that
$$ (P_1)_{lk}^{ij} = -a\hat q^{ij}q^{kl} (P_1)_{kl}^{ji} \not= 0.\eqno(7.1)
$$
\no This defines an elementary deformation if and only if the parameters
satisfy
the conditions
$$ q^{im} q^{jm} q^{mk} q^{ml} = \cases{a^x,&case 1, \cr
                                    a^{-x}, &case 2; \cr}\,\, \, x = \delta_m^i
- \delta_m^j, \, m = 1,2 \ldots N. \eqno(7.2)
$$
\no The ``exceptional series" of elementary, first order deformations exists
only if $a^3 = 1$. Let $i,j,k$ be neighbors in the natural numbers, with $i + 1
= j$. Let $P_1=0$  , except that either
$(P_1)_{kk}^{ij} = -aq^{ij}(P_1)_{kk}^{ji} \not= 0$, or else
$(P_1)_{ij}^{kk} = -q^{ij}(P_1)_{ji}^{kk} \not= 0$, but not both. This defines
an elementary deformation if and only if the parameters satisfy
$$
\eqalign {&(P_1)_{kk}^{ij} \not= 0:\,\,\,
         (q^{km})^2q^{mj}q^{mi} = a^x,\,\,\, x = \delta_{mi} - \delta_{mj};
           \cr
 &(P_1)_{ij}^{kk} \not= 0: \,\,\,
         (q^{km})^2q^{mj}q^{mi} = a^x,\,\,\,  x = \pm (\delta_{mk} -
\delta_{mi}). \cr} \eqno(7.3)
$$
\no The two signs in the last line apply when $k = i-1, k = j+1$, respectively.
There are no other first order, elementary deformations.  The elementary
deformations are formal and exact.

\vfill
\eject

\no {\stepone 8. The Classical Limit}
\b All the results obtained  here,  for twisted quantum $gl(N)$,  have direct
application to twisted, quantum $sl(N)$.  The connection between these two was
explained by Schirrmacher in [5] and is discussed also in [10].  In this
section
we shall take the classical limit and confront our results for $gl(N)$  with
the
classification, by Belavin and Drinfeld [1], of the classical r-matrices for
$sl(N)$. Strictly, this is possible only under additional conditions  on the
parameters, namely
$$
\prod_iq^{ij}a^j = a^{(N+1)/2};\eqno(8.1)
$$
\no we therefore assume that these relations hold,  although they do not
interfere directly with the following calculations.

The deformed quantum P-algebras of the principal series are semiclassical. The
classical r-matrix is defined by expanding the parameters
$$ a = 1 + h, \, q^{ij} = 1 + hp^{ij}, \, i<j, \eqno(8.2)
$$
\no and the R-matrix,
$$ R_{lk}^{ij}:=P_{kl}^{ij} + \epsilon (P_1)_{kl}^{ij}, \eqno(8.3)
$$
\no in powers of $h$,
$$ R = 1 - hr_\epsilon + O(h^2), \, r_\epsilon = r + \epsilon \delta r.
  \eqno(8.4)
$$
\no Here $r$ is the r-matrix for twisted (= multiparameter) $gl(N)$,
$$ r = \sum_{i<j} M_j^i \otimes M_i^j + r_0, \eqno(8.5)$$
\noindent with
$$r_0 := \sum_{i<j} \bigl(p^{ij} M_j^j \otimes M_i^i - (1 + p^{ij})
  M_i^i \otimes M_j^j \bigr), \eqno(8.6)
$$
\no and the perturbation is
$$ h\delta r_{kl}^{ij} = (P_1)_{lk}^{ij}, \,\,\, {\rm or} \, \,\,
  h\delta r = \sum_{i,j,k,l} (P_1)_{lk}^{ij} \,\,M_k^i \otimes
M_l^j.\eqno(8.7)
$$
\no Here $M_k^j$ is the matrix with the unit in row-$k$, column-$j$, all the
rest zero.

We examine the classical limit of an elementary, first-order deformation.  Fix
the notation as in Theorem 4, the expression for $\delta r$ is, up to a
constant,
$$
\delta r = M_k^i \otimes M_l^j - M_l^j \otimes M_k^i .\eqno(8.8)
$$
\no The diagonal matrices $(M_i^i), \, i=1,\ldots,N$, will be taken as a basis
for a ``Cartan subalgebra of $gl(N)$".  The upper triangular matrices form the
subspace of positive roots and the matrices $M_i^j$ with $i-j= \pm 1$ are the
simple roots.  [We have abused the notation by extending the notion of roots
from $sl(n)$ to $gl(n)$ and by introducing both positive and negative ``simple"
roots.]  The conditions on the indices that are spelled out in Theorem 4 insure
that all the roots appearing in (8.8) are simple, and that a positive root is
paired with a negative root and {\it vice versa}.  Let us find out what is the
meaning of the restriction (7.2) on the parameters.

Let $\Gamma_1(\Gamma_2)$ be the root space spanned by $M_k^i(M_j^l)$ and
$\tau: \Gamma_1 \rightarrow \Gamma_2$ the mapping defined by $\tau(M_k^i) =
M_j^l$.  Consider the equation
$$ (\tau \alpha \otimes 1 + 1 \otimes \alpha) r_0 = 0, \eqno(8.9)
$$
\no where $\alpha = M_k^i$ and $(\alpha \otimes 1) H \otimes H' =
\alpha(H)H'$ and $ (1 \otimes \alpha)H \otimes H' = \alpha(H')H$,
$H$ and $H'$ in the diagonal subalgebra of $gl(N)$.  We find that this equation
is equivalent to
$$ p^{lm} + p^{km} + p^{mi} + p^{mj} =  \delta_m^j - \delta_m^i,\eqno(8.10)
$$
\no  This equation is precisely the first order analog of Eq. (7.2), while
(8.9)
is a condition (Eq.~(6.7)) of Belavin and Drinfeld [1], applied to the case of
the elementary deformation (8.8). We have thus established that the conditions
(7.2)  on the parameters lift the invariance condition  of ref. [1] to the
quantum algebra.

\ve
\no {\stepone 9. A Family of Exact Deformations.}
\b       The elementary deformations are exact.  Combinations of elementary
deformations are not always exact and corrections of order
$\epsilon^2$ are needed. The interesting question is whether, for any first
order deformation, the requirements on the next order correction has a solution
(and so on). The case
$N = 2, a = 1$, shows that this is not always the case. It seems likely that
the
question has a cohomological formulation.

The dual Lie algebra of a standard, simple Lie bialgebra,  defined by the
simplest classical r-matrix, is solvable.  Since, in the application to
solvable
field theories, this dual  algebra is the algebra of principal dynamical
variables,  the nature of the physical applications is strongly limited. The
twisted forms of quantum $gl(N)$ do not offer anything new in this respect, and
indeed, multiparameter quantum $gl(N)$ has been stigmatized as a mere ``gauge
transformation" of the standard version [7].  Esoteric quantum $gl(N)$ is quite
different.  We describe here a special case that may represent  the strongest
departure from a solvable dual.  The special cases $gl(3)$ and $gl(5)$ were
described in  some detail in a preprint [8], and $gl(2n-1)$ in a conference
report [9]. Here is ``esoteric quantum
$gl(2n-1)$" with some additional information.

        We begin with multiparameter $gl(2n-1)$, with the special  values of
the
parameters
$$
\matrix{q^{ij} = 1/q, &i<j, &i+j \not= 2n,  &(p^{ij} = -1/2)\hfill \cr q^{ij} =
1/ q^2, &i<j, & i+j = 2n, &(p^{ij} = -1)\hfill\cr}\eqno(9.1)
$$
\no and $a = q^2$, with the R-matrix
$$
\eqalign{R_0 = &{\displaystyle \sum^{2n-1}_{i=1}} M_i^i \otimes M_i^i +
{\displaystyle \sum_{i<j<2n}}(1-q^2)M_j^i \otimes M_i^j +  {\displaystyle
\sum_{i \neq j,i +
   j \neq 2n}}q M_i^i \otimes M_j^j  \cr &\hskip1in + \,\, {\displaystyle \sum
_{i<n}} M_i^i\otimes M_{i'}^{i'} +  {\displaystyle\sum _{i<n}} q^2 M^{i'}_{i'}
\otimes M^i_i, \cr}\eqno(9.2)
$$
\no The R-matrix for the deformation that  we have called esoteric quantum
$gl(2n-1)$ is
$$
\eqalign {R = &R_0 + R_1, \cr R_1 = &\displaystyle \sum_{k<i<n}(\mu'_iM_i^n
\otimes M_{i'}^n  + \mu_iM_{i'}^n \otimes M_i^n)\, \cr &+
\sum_{k<i<j<n}(\lambda'_{ij}M_i^j \otimes M_{i'}^{j'} + q^2 \lambda_{ij}
M_{i'}^{j'} \otimes M_i^j) \cr}\eqno(9.3)
$$
\no We have introduced the notation
$$  i' := 2n-i, \,\,\, {\rm for} \,\,\, 0<i<n.\eqno(9.4)
$$ The relations of the corresponding quantum (anti-) plane  are the same as in
(2.3, 2.4) when $i+j \not= 2n$  and, for $j<n$:
$$
\eqalign {[x^j,x^{j'}]_{1/q^2}
 &= q^{-2} \sum_{i<j} \lambda_{ij}x^{i'}x^i
 = -\sum_{i<j} \lambda'_{ij}x^ix^{i'}, \cr   [\theta^j,\theta^{j'}]
 &= - \sum_{i<j} \lambda_{ij}\theta^{i'}\theta^i
 = - \sum_{i<j} \lambda'_{ij}\theta^i\theta^{i'}, \cr
\theta^n \theta^n &= - \sum_{i,n} \mu'_i \theta^{i'}\theta^i, \cr [x^j,
\theta^{j'}]_{1/q^2} &=q^{-2} \sum_{i<j} \lambda_{ij}\theta^{i'}x^i
 + (q^{-2} - 1) \theta^jx^{j'}, \cr [\theta^j,x^{j'}] &= - q^{-2} \sum_{i<j}
\lambda'_{ij}\theta^ix^{i'}, \cr [\theta^n,x^{j'}] &= -\sum_{i<j} \lambda'_{ij}
\theta^ix^{i'}, \cr [\theta^n,x^n]_{q^2} &= -q^2 \sum_{i<j} \bigl(\mu'_i
\theta^ix^{i'}  + \mu_i \theta^{i'}x^i  \bigr). \cr }\eqno(9.5)
$$
\no Here $[A,B]_q := AB -qBA.$ The braid relation is easily imposed by
invoking
Theorem 3, with the result that  the coefficients $\mu, \mu', \lambda $ and
$\lambda'$  must satisfy the conditions
$$
\matrix{
\mu'_i &=& - q^{2(i-n)} \mu_i,\hfill&i<n,\hfill \cr
\mu_j \lambda_{ij} &=& (1-q^{2})q^{2(i-j)} \mu_i,\hfill& i<j<n,\hfill \cr %
\mu_j \lambda'_{ij} &=& (q^2 - 1)\mu_i,\hfill&  i<j<n.\hfill\cr}
\eqno(9.6)
$$
\no The parameters $(\mu_i), i = 1,2,...,n, $ remain arbitrary except  that
$\mu_i = 0$ implies $\mu_{i-1} = 0,$

Esoteric quantum $gl(2n-1)$ (in the Woronowicz picture) is  generated by the
matrix elements $z_i^j$ of a matrix $T$,  of dimension $2n-1$, with relations
given by
$$  P(T\otimes T) = (T\otimes T)P.\eqno(9.7)
$$
\no These relations are complicated and we refer to [9].

\vfill
\eject

{\bf Acknowledgements.}

This work was supported in part by the National Science Foundation. C.F. thanks
the Research Committee of the Academic Senate, UCLA, for support. A.G. thanks
DGCICYT (Spain) for support.

\vskip.5in

\no {\bf References}
\b
\frenchspacing
\item{1.} A.A. Belavin and V.G. Drinfeld,  Sov. Sci. Rev. Math. {\bf 4} (1984),
93-165.
\item{2.} V.G. Drinfeld, {\sl Quantum Groups}, Intern. Congr. Math. Berkeley
1986, 798-820.
\item{3.} M. Jimbo,  Commun. Math. Phys. {\bf 102} (1986), 537-547.
\item{4.} P. Truini and V.S. Varadarajan,  Lett. Math. Phys. {\bf 21} (1991),
287;  P. Bonneau, M. Flato and G. Pinczon, Lett. Math. Phys.  {\bf 25} (1992),
75-84.
\item{5.} A. Sudbery, J. Phys. A{\bf 23} (1990), L697;  A. Schirrmacher, Z.
Phys. C{\bf 50} (1991) 321;   N.Y. Reshetikhin, Lett. Math. Phys. {\bf 20}
(1990), 331-336.
\item{6.} L. Woronowicz, Publ. Res. Ins. Math. Sci.  Kyoto Univ., {\bf 23}
(1987) 117-181.
\item{7.} N.Y. Reshetikhin, ref. 5.
\item{8.} C. Fronsdal and A. Galindo, {\sl  New Quantum Groups,  Quantum Planes
and
   Applications}, preprint UCLA/92/TEP/ and FT/UCM/4/92.

\item{9.} C. Fronsdal and A. Galindo, {\sl The Universal $T$-Matrix},
Proceedings of the 1992 Joint Summer Research Conference  on Conformal Field
Theory, Topological Field Theory and Quantum Groups,   Holyoke, June 1992, to
be
published (preprint UCLA/93/ TEP/2);  Lett. Math. Phys. {\bf 27} (1993),
59-71.
\item{10.} C. Fronsdal, {\sl Universal $T$-matrix  for Twisted Quantum
$gl(N)$},
Proceedings of the NATO Workshop, San Antonio  February 1993.
\item{11.} E. Cremmer and J.-L. Gervais, Commun.Math.Phys.  {\bf 134} (1990)
619-632.
\item{12.} C.Fronsdal, {\sl Cohomology and quantum groups}, Proceedings  of the
XXX'th Karpacz Winter School of Theoretical Physics, February 1994, to be
published.
\bye